# A Line Search Algorithm for Multiphysics Problems with Fracture Deformation

Author: Ivar Stefansson, Center for Modeling of Coupled Subsurface Dynamics, Department of Mathematics, University of Bergen.

**Abstract**

Models for multiphysics problems often contain strong nonlinearities. Including fracture contact mechanics introduces discontinuities at the transition between open and closed or sliding and sticking fractures. The resulting system of equations is highly challenging to solve. The naïve choice of Newton's method frequently fails to converge, calling for more refined solution techniques such as line search methods.

When dealing with strong nonlinearities and discontinuities, a global line search based on the magnitude of the residual of all equations is at best costly to evaluate and at worst fails to converge. We therefore suggest a cheap and reliable approach tailored to the discontinuities. Utilising adaptive variable scaling, the algorithm uses a line search to identify the transition between contact states. Then, a solution update weight is chosen to ensure that no fracture cells move too far beyond the transition.

We demonstrate the algorithm on a series of test cases for poromechanics and thermoporomechanics in fractured porous media. We consider both single- and multifracture cases and study the importance of proper scaling of variables and equations.

# 1. Introduction

This paper concerns solution strategies for numerically solving strongly nonlinear and non-smooth equation systems. The primary motivation is multiphysics problems involving fracture contact mechanics in porous media. The developed methods may, however, be relevant for other problems involving friction and contact constraints.

The Newton-Raphson method is the go-to method for solving systems of nonlinear equations due to its quadratic convergence. However, convergence is only guaranteed locally, which is especially restrictive in the non-smooth setting. This motivates using globalisation techniques ensuring convergence. The pursuit of appropriate solution strategies is an emerging and topical research area (White et al., 2019).

Globalisation schemes are well studied in the field of optimisation, including application to problems originating in PDEs. Following Nocedal and Wright (1999), we distinguish between two main families of methods. Trust region methods prescribe a local region around the current iterate in which the update is sought, e.g. using the descent direction provided by Newton's method. By contrast, line search methods first compute the update direction, before searching along that direction for the best solution according to some metric, e.g. minimising the residual. They therefore tend to be less intrusive into the solution algorithm, and are the path pursued in this paper.

Inequality conditions such as those arising in contact mechanics lead to *constrained* optimisation problems, requiring modification of the globalisation scheme. A conceptually straightforward approach is to replace the objective function representing the system of equations by a merit function additionally incorporating information about the constraints. This is achieved using additional variables penalising violation of the constraints. As discussed by Hiermeier (2020), the penalty approach comes with the inherent challenge of choosing and adapting the penalty parameters which balance the two goals of minimising the objective function and honouring the constraints.

We expect the non-smooth constraints to be the main source of difficulty for the Newton solver. Indeed, in their review of numerical solution of contact problems, Acary et al. (2018) point out that the standard methods may need additional criteria to terminate the line search. Moreover, in multiphysics problems, valuable information contained in the constraint functions may be obscured by the rest of the residual.

Together with the cost of performing a line search on the global residual, these considerations motivate a more targeted way of searching along the Newton direction. Herein, we pursue a line search exploiting knowledge of the problem's irregularities. The search is based on the discontinuities of the contact mechanics relations and is evaluated for the fracture cells only. It prevents updates from passing far beyond a singular point of the constraints, as suggested in the context of multiphase flow by Khebzegga et al. (2021) and Moyner (2017). We combine this line search with adaptive scaling of the constraint conditions to achieve an algorithm which is both efficient and dependable.

The rest of the paper is structured as follows: Section 2 outlines the mathematical model for frictional contact mechanics and mixed-dimensional thermoporomechanics, highlighting prominent nonlinearities and non-smoothness. The new solution strategy is described in Section

3 and illustrated and tested by the simulation results of Section 4. Finally, we offer concluding remarks in Section 5.

## 2. Mathematical Model

We consider a mixed-dimensional discrete fracture-matrix model as described in (Boon et al., 2021), where the domain is partitioned into subdomains $\Omega_i$ of dimension $d_i$. Each pair of neighbouring subdomains separated by one dimension is connected through an interface $\Gamma_j$. We detail only the equations of particular relevance to the suggested algorithm herein, referring to the paper by Stefansson et al. (2024) for the full mixed-dimensional models: The poromechanical model consists of Eqs. (1), (4)–(10), (32) and (34), ignoring the temperature terms. We obtain the thermoporomechanical by adding Eqs. (2) and (33). We close the models by including the relevant boundary conditions and constitutive laws as defined in the paper's Sections 3.4 and 3.5, replacing Eq. (28) by Eq. (*2*) defined below.

The fracture contact mechanics being of particular importance in the present context, we repeat the central expressions here. Following (Hüeber et al., 2008), we formulate the fracture deformation constraints as non-smooth complementarity functions:

$$\begin{aligned} C_\perp(\sigma, [\![u]\!]) &= -\sigma_\perp - \max\{0, -\sigma_\perp - c([\![u]\!]_\perp - g)\}, \\ C_\parallel(\sigma, [\![u]\!]) &= \chi_o \sigma_\parallel - (1 - \chi_0)\left[\sigma_\parallel \max\{b, \|\sigma_\parallel + c[\![\dot u]\!]_\parallel\|\} - b(\sigma_\parallel + c[\![\dot u]\!]_\parallel)\right]. \end{aligned} \quad (1)$$

Here, $\sigma$ and $[\![u]\!]$ denote contact traction and displacement jump between the two sides of the fracture, with the dot indicating increment between successive (time) states. $b(\sigma) \coloneqq -F\sigma_\perp$ is the friction bound with friction coefficient $F$, herein assumed to be constant. The characteristic function $\chi_o \coloneqq b \leq 0$ indicates whether the fracture is open, while $c$ is a numerical parameter. Subscripts $\perp$ and $\parallel$ denote the normal and tangential direction of the fracture, respectively. Finally, $g$ is the gap function, i.e., the distance between the fracture surfaces when in mechanical contact. $g$ accounts for shear dilation according to

$$g = \tan\phi \, \|[\![u]\!]_\parallel\|. \quad (2)$$

The dilation angle $\phi$ determines the strength of the coupling between tangential and normal displacement and by extension the effect on the hydraulic aperture $a \coloneqq a_{res} + [\![u]\!]_\perp$. To emphasize the strength of the nonlinear coupling, we specify that the fracture permeability is related to $a$ by the cubic law.

## 3. Solution strategy

The spatial discretisation uses multi-point finite volume methods for stress and diffusive fluxes (Nordbotten & Keilegavlen, 2021) and a first-order upwind scheme for advective fluxes. To isolate the effect of the contact mechanics solution strategy, we use the TPFA scheme in the fractures. Thus, the nonlinearity in the fracture permeability may be treated fully implicitly as described in (Stefansson & Keilegavlen, 2023).

The implementation is provided in the PorePy simulation toolbox for fractured porous media, which is described by Keilegavlen et al. (2021) and Stefansson et al. (2024), which also contains a description of the parts of the solution strategy and discretization not detailed herein.

Denoting the solution vector by $x \in \mathbb{R}^n$, we write the system of discretised nonlinear equations by

$$r(x) = 0. \tag{3}$$

The strongly nonlinear and tightly coupled nature of $r$ causes significant difficulties in solving the system. The following sections presents an algorithm centred around the fracture deformation equations. It consists of an adaptive scaling and a line search local to the fractures.

## 3.1. Scaling

As for most solution algorithms, the present one relies on judgements on the magnitude of (updates) of the variables and equations. Wishing to avoid the error-prone approach of combinations of relative and absolute tolerances throughout the algorithm, we attempt to scale $\sigma$ and $[\![u]\!]$ towards unity.

An educated a priori guess about the expected magnitude of the unscaled variables can be based on the driving forces (boundary conditions, source and sink terms). The simulations in Section 4 all have a Dirichlet boundary condition for displacement as a primary driving force. In this case, we can define the characteristic displacement $u_c$ to equal the boundary value. Furthermore, since the examples are defined on unit cube domains, we set the characteristic traction to $\sigma_c = Eu_c/1$, with $E$ being Young's modulus in the matrix. Note that in case of elastic normal deformation of the fracture as described by e.g. (Bandis et al., 1983), the fracture's normal stiffness may be a more appropriate choice. Given a choice of $\sigma_c$ and $u_c$, we can replace $\sigma$ in Eq. (1) by the scaled traction variable $\tilde{\sigma} = \sigma/\sigma_c$. Instead of using a scaled displacement, we choose $c = 1/u_c$, thus scaling $g$ as well as $\boldsymbol{u}$.

We stress that, however experienced the practitioner, this guess will contain high uncertainty, typically several orders of magnitude, for complex, real world applications. Therefore, we suggest an adaptive scaling to be used in combination with the line search as described in the subsequent sections. We first define a cell-wise scaling estimate $s^\nu$. Aiming to emphasize high values without ignoring lower ones, we then employ the p-mean with $p = 5$ to the local estimate for all fracture cells:

$$s^\nu = \|\tilde{\sigma}^\nu\| + \|c([\![\boldsymbol{u}]\!]^\nu - \boldsymbol{n}g^\nu)\|, \tag{4}$$
$$s = \left(\frac{\sum (s^\nu)^p}{\#n_f}\right)^{1/p}.$$

Here, $\nu$ denotes individual cells, $\#n_f$ the total number of fracture cells and $\boldsymbol{n}$ is the normal vector of the fracture. The scale $s$ is fixed to its value at the previous iteration. For numerical robustness, we recommend capping $s$, herein between $10^{-8}$ and $10^8$.

## 3.2. Line search

To pose the problem in optimisation terminology, we define an objective function $f$, which is usually taken as some norm of $r$. Closely following Nocedal and Wright (1999), we consider the standard choice

$$f = \frac{1}{2}\|r\|^2$$

and seek the minimizer $x^*$ satisfying

$$\nabla f(x^*) = 0.$$

Newton's method defines the following iteration scheme to obtain the Newton update at iteration $k$:

$$p^k = -J^{-1}(x^k)r(x^k),$$

with $J(x)$ denoting the Jacobian of $r$. A line search method now involves determining a beneficial step length $0 < \alpha_t^k \leq 1$ and updating the iterate according to

$$x^{k+1} = x^k + \alpha_t^k p^k.$$

### 3.2.1. Global, residual-based line search

The basic line search algorithms are based on minimizing $f$ along $p^k$ subject to some conditions ensuring convergence. These approaches inadvertently add significantly to the computational cost of each iteration. This is especially true for non-smooth problems, which may require dense sampling of trial values of $f$. However, we include a simple residual-based algorithm as a comparison to the approach described in the subsequent section. To somewhat alleviate computational cost, we sample at a limited number of points and interpolate using the monotone cubic spline scheme proposed by Fritsch and Carlson (1980) as suggested by Moyner (2017) in a similar context involving discontinuities.

### 3.2.2. Local, constraint-based line search

To design a more targeted line search, we borrow concepts used in the context of transport problems by Moyner (2017) and (Pour et al., 2023). Since the discontinuities are due to the maximum functions, we design an algorithm based on the relative values of their arguments: If the computed Newton update leads to a reversal of maximum arguments, we seek a weight corresponding to a damped update just beyond the transition point.

To facilitate efficient weight computation, we introduce state indicator functions related to $\max(\phi, \varphi)$ which are linear in the arguments $\phi$ and $\varphi$ and change sign at the discontinuity. For comparison, we define both a constantly and adaptively scaled version:

$$i_c = (\theta - \varphi) \tag{5}$$
$$i_a := \frac{i_c}{s} = (\theta - \varphi)/s.$$

Note that the magnitude range of the arguments should be similar, to ensure robustness of the root seeking. Scaling $i$ to approximate unity implies the tolerance introduced below is independent of the problem specifics. For the normal constraint of Eq. (*1*), we compose $i_{c\perp}$ from

$$\theta_\perp = -\tilde{\sigma}_\perp, \tag{6}$$
$$\varphi_\perp = c(u_\perp - g),$$

while in the tangential case, we have

$$\theta_\parallel = |\tilde{\sigma}_\parallel - c_\parallel \dot{u}_\parallel|, \tag{7}$$
$$\varphi_\parallel = b(\tilde{\sigma}).$$

To avoid restricting open cells by the tangential constraint, we multiply by the Heaviside function $H(i_{c\perp}) = i_\perp > 0$, obtaining $i_{c\parallel} = (\theta_\parallel - \varphi_\parallel)H(i_{c\perp})$.

Denoting a trial weight by $\alpha_t^k$, we define the transition indicator as

$$t(\alpha_t^k) := -\text{sgn}\{i(x^k) \cdot i(x^k + \alpha_t^k p^k)\} \cdot |i(x^k + \alpha_t^k p^k)|. \tag{8}$$

The first factor identifies cells transitioning between contact states, while the second factor limits the degree to which they transition into the new state. Introducing a constraint violation tolerance $\delta > 0$, we compute cell-wise weights $\alpha_t^\nu$ for the cells satisfying $t^\nu > \delta$ such that

$$i(x^k + \alpha_c^{\nu,k} p^k) + \delta \text{sgn}\left(i(x^k)\right) = 0, \tag{9}$$

again using the interpolation line search described above. The global trial weight $\alpha_t$ is taken as the minimum among all $\alpha_t^\nu$.

The tolerance $\delta$ may allow multiple fracture cells to transition within one iteration. While significantly speeding up convergence, this can in rare cases lead to loss of convergence if too many cells transition. Therefore, we recursively impose a tightening of the tolerance whenever a fracture's number of transitioning cells $\#t_i = \sum_{\nu \in \Omega_i} t^\nu > 0$ is high relative to the fracture's number of cells, $\#\nu_i$, i.e.,

$$\#t_i > \max(1, \gamma \cdot \#\nu_i), \tag{10}$$

with the cutoff 1 intended for poorly resolved meshes and $\gamma$ denoting the relative tolerance. That is, we halve $\delta$ and recompute $\alpha_t^\nu$ according to Eq. (9) as long as Eq. (*10*) is satisfied.

# 4. Simulations

This section demonstrates the algorithm's efficiency and reliability for a range of test cases, which differ along five dimensions as detailed in the subsections. Furthermore, we compare the suggested constraint based line search with adaptive scaling (CLS $i_a$) to three approaches: no line search (No LS), a residual-based line search (RLS) and a constraint line search with constant scaling using $i_c$ from Eq. (5) (CLS $i_c$).

The tolerances of the constraint-based algorithm are set to $\delta = 0.3$ and $\gamma = 0.2$ throughout. We run both simulations for a single time step of length $10^6$ s.

We compare the results in terms of number of nonlinear iterations, marking runs which either did not converge within one hundred iterations or diverged (i.e., the residual containing infinite entries or similar) in grey and labelling them with NC or Div, respectively. While we do not report run times due to a non-optimised implementation, we note that, unlike the local line search, the global line search adds perceptively to overall computational cost.

## 4.1. Single fracture

In the first suite of test cases, we include both poromechanics and thermoporomechanics, thus demonstrating applicability to different physics. We vary the dilation angle $\phi \in \{0.1, 0.2\}$, which influences the coupling strength from deformation to fracture flow, cf. Eq. (*2*). We prescribe two different mesh sizes $h \in \{1/6, 1/12\}$, thus testing efficiency with respect to the number of fracture cells. Finally, we vary the characteristic displacement scaling employed in Eq. (*1*) to investigate the robustness in cases where this quantity is difficult to assess a priori. Since using differently scaled variables affects the residual magnitude, we employ a convergence criterion based on the normalised $L^2$ norm of the nonlinear increment, $\|p^k\|/\sqrt{n} < 10^{-10}$, thus allowing comparison across cases.

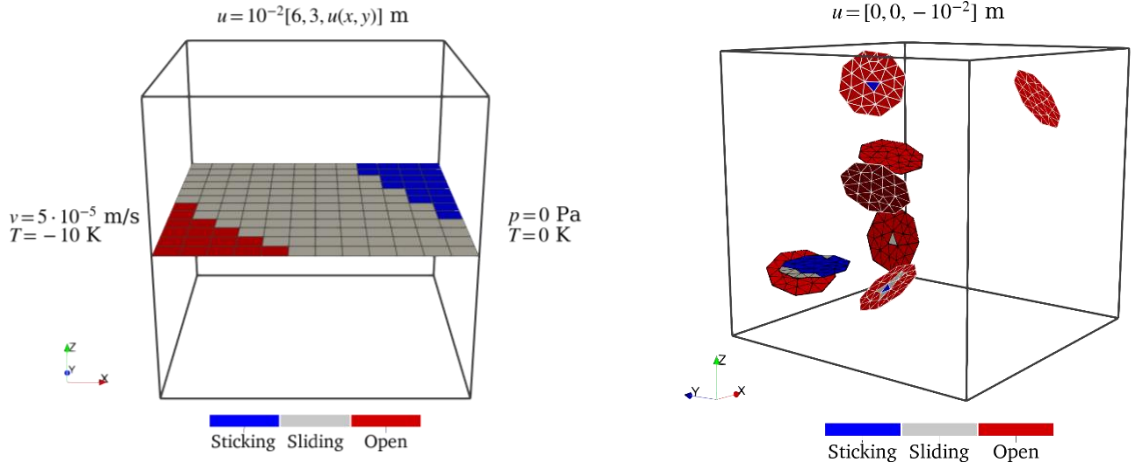

*Figure 1: Geometry, selected boundary conditions and fracture deformation solution for the poromechanics $\phi = 0.1$ cases. On the left, we show a Section 4.1 case, with $v$ denoting specific fluid discharge. In the Section 4.2 case on the right, fractures included in the four-fracture case are shown with black edge lines and those unique to the eight-fracture case with white lines. The mostly hidden fracture is fully open.*

This example is defined on a unit cube domain with a single throughgoing fracture as shown in *Figure 1*. There is inflow on the left and outflow on the right fracture boundary, and no-flow conditions elsewhere. We prescribe zero and heterogeneous displacement values on the bottom and top boundary, respectively, resulting in fracture deformation containing both sticking, sliding and open fracture cells as shown in *Figure 1*. The remaining parameters are listed in *Table 1*. While these do not represent any particular physical setting, some characteristics, such as relatively low values for permeability and high ones for stiffness, contribute to retaining the relevant relative importance of the different terms and couplings in the equation system.

Taking the undamped method (No LS) as the baseline, the results reported in Figure 2 show that the thermoporomechanical problem is more challenging than the poromechanical. Similarly, increasing the number of cells or the dilation angle adds to the difficulty. Both the residual-based and the constantly scaled constraint-based line search (RLS and CLS $i_c$) obtain convergence in some of the cases where the reference method does not. However, neither is robust with respect to the characteristic displacement scaling, with increased iteration counts or failure to converge in several cases. In contrast, the adaptively scaled constraint-based method (CLS $i_a$) reliably converges with an iteration count which is constant with respect to the characteristic displacement. The number of iterations is also consistently low and competitive with all other converging algorithms. Moreover, the iteration count hardly increases with the number of fracture cells, which equals $1/h^2$.

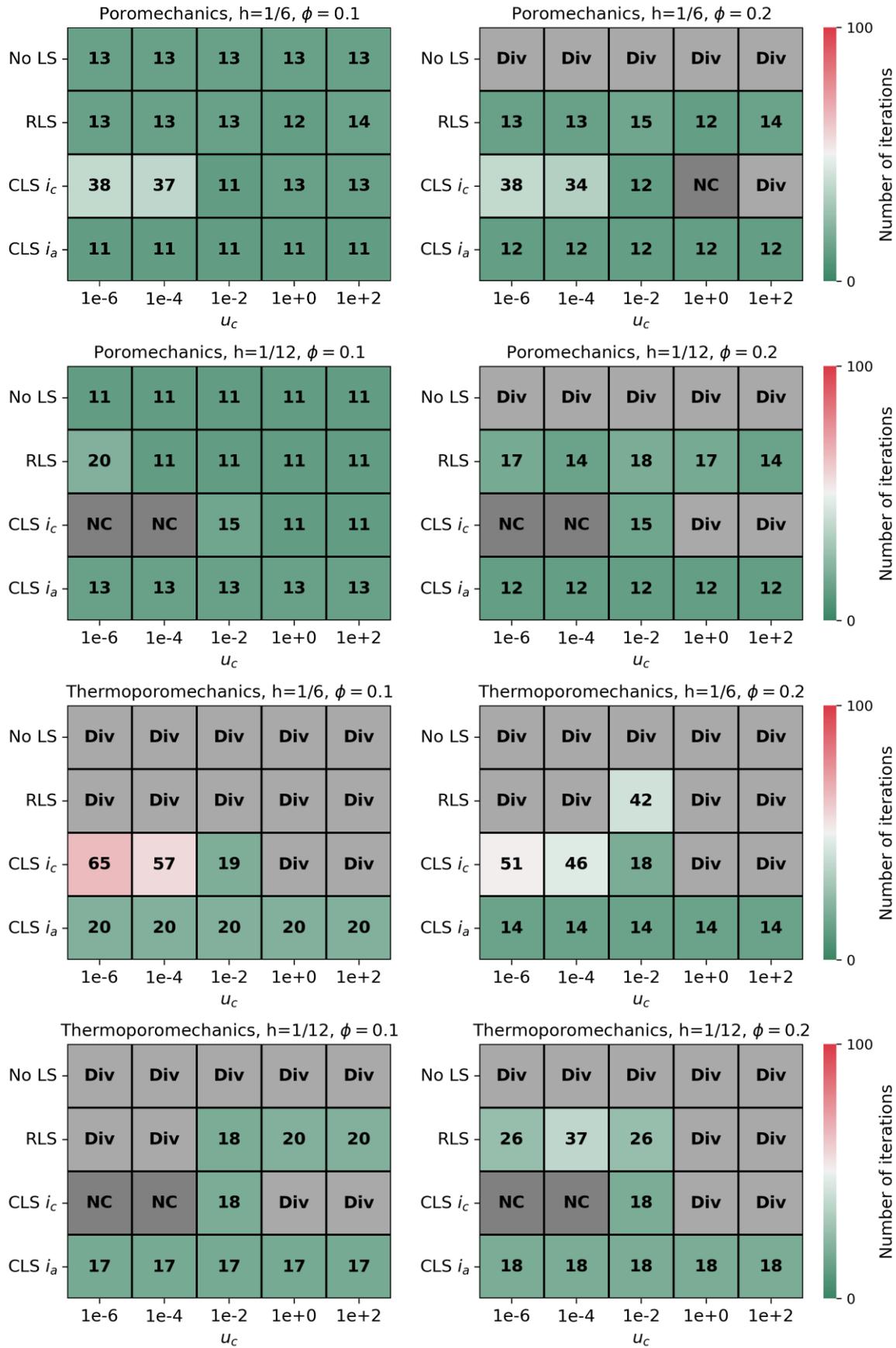

Figure 2: Iteration counts for the examples of Section 4.1. NC and Div indicate simulations which did not converge within one hundred iterations or diverged.

*Table 1: Simulation parameters.*

| Fluid parameters | | |
|---|---|---|
| Compressibility | $1.0 \cdot 10^{-6}$ | 1/Pa |
| Density | 1.0 | kg/m$^3$ |
| Normal thermal conductivity | 1.0 | W/m/K |
| Reference pressure | 0.0 | Pa |
| Specific heat capacity | 100.0 | J/kg/K |
| Reference temperature | 0.0 | K |
| Thermal conductivity | 1.0 | W/m/K |
| Thermal expansion | 0.01 | 1/K |
| Viscosity | 0.1 | Pa s |
| Inlet/outlet pressure | $1.5 \cdot 10^5 / -1.0 \cdot 10^5$ | Pa |
| Inlet/outlet temperature | $-10.0/0.0$ | K |
| Solid parameters | | |
| Biot coefficient | 0.8 | - |
| Density | 1.0 | kg/m$^3$ |
| Dilation angle | 0.1, 0.2 | - |
| Friction coefficient | 1.0 | - |
| First Lamé parameter | $2.0 \cdot 10^6$ | Pa |
| Normal permeability | $1.0 \cdot 10^{-6}$ | m^2 |
| Permeability | $1.0 \cdot 10^{-8}$ | m^2 |
| Porosity | 0.01 | - |
| Residual aperture | $1.0 \cdot 10^{-3}$ | m |
| Shear modulus | $2.0 \cdot 10^6$ | Pa |
| Specific heat capacity | 100.0 | J/kg/K |
| Reference temperature | 0.0 | K |
| Thermal conductivity | 1.0 | W/m/K |
| Thermal expansion | $1.0 \cdot 10^{-3}$ | 1/K |

## 4.2. Multiple fractures

The second suite also considers two physical models and $\phi \in \{0.1, 0.2\}$, but fixes $u_c = 0.01$. Additionally, we run with both four and eight fully immersed, randomly oriented fractures in the unit cube domain (see supplementary material for their vertex coordinates). We prescribe zero fluid and energy flow at all external boundaries. For the momentum balance, we use zero and a constant compressive displacement at the bottom and top, respectively, and zero traction elsewhere. In the centremost cell of each fracture, we prescribe a constant pressure and temperature value, which we pick from two values (which may be interpreted as mimicking injection and production wells). This results in deformation regimes for individual fractures ranging from mostly open through sliding to sticking, depending on fracture orientation and assigned pressure. The remaining parameters are as in the previous section. Figure 3 shows the fracture geometry, as well as illustrating the contact mechanical state of the poromechanical simulation with $\phi = 0.1$. Since the equations are scaled consistently for all cases, we use a residual-based convergence criterion $\|r^k\|/\sqrt{n} < 10^{-10}$.

The iteration counts shown in Figure 3 demonstrate similar trends in convergence behaviour with respect to physical models and $\phi$. Unsurprisingly, increasing the number of fractures also reduces the likelihood of convergence. Again, the residual-based search offers significant but unsatisfactory improvement over the undamped Newton algorithm. The two constraint-based methods consistently converge with similar iteration count, which indicates a quite accurate estimate of $u_c$. The number of iterations scales very modestly with the number of fractures. This indicates that employing a permissive tolerance $\delta$ allows multiple fractures to partly transition within the same iteration.

We emphasize that no convergence issues arose with the relatively lenient tolerances in any test case. Since reducing the tolerances will at some point increase the number of iterations, we recommend this more aggressive choice.

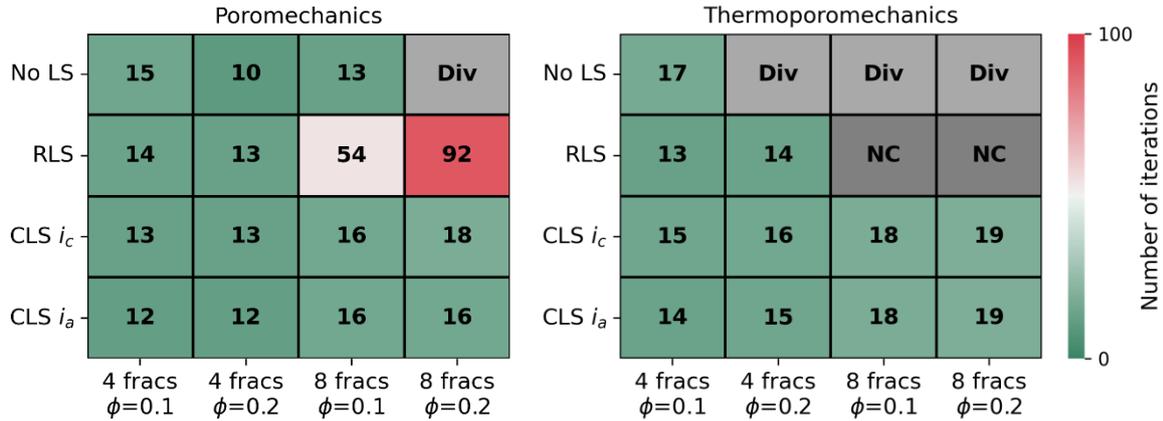

*Figure 3: Iteration counts for the Section 4.2 test case with multiple fractures. NC and Div indicate simulations which did not converge within one hundred iterations or diverged.*

# 5. Conclusion

We have presented a line search algorithm for solving multiphysics problems involving fracture contact mechanics. The nonlinear and non-smooth nature of the equations is addressed in a targeted manner producing a simple, efficient and reliable algorithm. The algorithm consists of adaptive variable scaling and a line search based on the non-smooth part of the fracture contact mechanics equations.

Numerical simulations demonstrate applicability to both poromechanical and thermoporomechanical problems and ability to deal with variables of challenging scales. The suggested approach consistently converges, as opposed to both the standard Newton method and less targeted line search approaches.

We also assess efficiency in terms of number of nonlinear iterations. The results show competitiveness with the alternative approaches for the cases where the latter converge, indicating that the line search hardly introduces any reduction in convergence rate. Moreover, we obtain very favourable scaling with respect to both number of fracture cells and the number of fractures.

# Acknowledgements

This work was funded by the VISTA program, The Norwegian Academy of Science and Letters and Equinor.

## Availability of Data and Code

The simulations were run using PorePy (Keilegavlen et al., 2021) and the run scripts are provided in a GitHub repository at https://github.com/IvarStefansson/A-Line-Search-Algorithm-for-Multiphysics-Problems-with-Fracture-Deformation/. The repository contains a setup for running with a development container based on Visual Studio Code for convenient reproduction code inspection. The full environment is available at Zenodo as a Docker image (Stefansson, 2024).